# A Lattice of Semigroup Inclusion Classes consisting of Unions of Varieties of Generalised Inflations of Rectangular Bands

## R. A. R. Monzo[1]

**Abstract**  We determine the lattice of inclusion classes under $[xyz \in \{xywz, xqyz\}; xy = xyxy]$.
.

**Keywords**  Inflation, Generalised inflation, Lattice, Inclusion class

## 1. Introduction

Following the review of T. Evans in 1971 the lattice of semigroup varieties has received much attention [2]. The subsequent survey of Shevrin, Vernikov and Volkov (2009) summarises some of the significant results and open questions [10]. Certain sub-lattices of the lattice of varieties of semigroups have been determined, such as the sub-lattice of varieties of bands and the sub-lattice of semigroups whose cube is a rectangular band and whose square is an inflation of its cube [1,3,4,5,8].

Ljapin noted that the notion of a semigroup variety can be extended to the notion of an "identical inclusion variety" [7]. Identical inclusion varieties were called "inclusion classes" in [9]. It will be clear from the definition given below that a semigroup variety is a semigroup inclusion class.

So one can extend the concept of the lattice $L$ of semigroup varieties to that of the lattice $L^*$ of inclusion classes of semigroups. Kalicki and Scott [6] determined the elements covering $\{0\}$, called atoms, in the lattice $L$. The atoms of $L^*$ were determined by Ljapin [7].

The sub-varieties of the variety of generalized inflations of rectangular bands were determined in [8]. Two of these sub-varieties are $[xyz = xywz ; xy = xyxy]$ and $[xyz = xwyz ; xy = xyxy]$. The union $I$ of these two varieties is not a variety. However it is a semigroup inclusion class. In fact we will prove that $I = [xyz \in \{xywz, xqyz\} ; xy = xyxy]$. The main result of this paper is the determination of the lattice of semigroup inclusion classes under $I$. We prove that every inclusion class under $I$ is a union of subvarieties of the variety of generalized inflations of rectangular bands. We also prove that, in general, not every inclusion class is a union of semigroup varieties.

## 2. Preliminaries

Throughout this paper $C = [L_k \subseteq M_k]_{k=1}^{t}$ will denote the collection of all semigroups S that satisfy the set of inclusions $\{L_k \subseteq M_k\}_{k=1}^{t}$, where each $L_k$ and $M_k$ is a finite set of words over some alphabet $\{x_\alpha\}_{\alpha \in \Omega}$ of variables indexed by $\Omega$. That is, for every substitution $* = (x_\alpha \to a_\alpha)_{\alpha \in \Omega}$ of every variable $x_\alpha$ by the element $a_\alpha \in S$, $L_k^* \subseteq M_k^*$ ($k \in \{1,2,\ldots,t\}$), where $L_k^* \subseteq M_k^*$ is the inclusion obtained from $L_k \subseteq M_k$ by making the substitution $*$.

If $C$ is an inclusion class then for every substitution $*$, for every word $L \in L_k$ there exists a word $W_{L,*} \in M_k$ such that $L^* = W_{L,*}^*$. Note that the word $W_{L,*}$ depends on L and $*$. When it is clear which substitution $*$ we are dealing with we may write $W_L$ instead of $W_{L,*}$. We will also use the symbol $\bullet$ to represent a particular substitution. If $a \in S$ then $*^a$ will denote the substitution that replaces every variable with the element a. In [9] the collection $C = [L_k \subseteq M_k]_{k=1}^{t}$ was called a ***semigroup inclusion class***. In [7] it was called an *identical inclusion variety of semigroups.* Clearly, as noted in [7], every semigroup variety $V$ is a semigroup inclusion class. We shall denote the collection of all inflations of semigroups in the variety $V$ and the collection of all generalized inflations of semigroups in $V$ as $\mathbf{I}V$ and $\mathbf{G}V$ respectively [9].

We also use the following notation:

$$Z = [xy = zw], L = [xy = x], R = [xy = y], RB = L \times R$$

$$\mathbf{GRB}_l = [xyz = xywz\,;\, xy = xyxy], B = [x = xx], \mathbf{GRB}_r = [xyz = xwyz\,;\, xy = xyxy]$$

Results 1, 2 and 3 below are well-known.

**Result 1.** $Z = [xyz = uv], L = [xyz = x], R = [xyz = z], \mathbf{I}L = [xyz = xw], \mathbf{I}R = [xyz = wz],$
$\mathbf{I}B = [xy = x^2 y^2], \mathbf{G}B = [xy = xyxy], \mathbf{G}L = [xyz = xy] = \{S: S^2 \in L\} = [xyx = xyw\,;\, xy = xyxy],$
$\mathbf{G}R = [xyz = yz] = \{S: S^2 \in R\} = [xyx = xyw\,;\, xy = xyxy], RB = [xyx = x] = [xyz = xz\,;\, x = xx]$
$\mathbf{I}RB = [xyz = xz], \mathbf{G}RB = [xy = xyzxy] = \{S: S^2 \in RB\}.$

**Definition.** If $x \in S \in RB$ then $x_1$ and $x_2$ will denote the first and second components of $x$ respectively. (Note that for any $\{x, y\} \subseteq S \in RB$, $x = (x_1, x_2)$ and $xy = (x_1, y_2)$.)

**Result 2.** $S \in [xy \in \{x, y\}]$ if and only if $S$ is a chain $C$ of disjoint semigroups $S_\alpha \in L \cup R (\alpha \in C)$ such that if $\alpha < \beta$, $x \in S_\alpha$ and $y \in S_\beta$ then $xy = yx = x$.

**Result 3.** If $S \in RB$ and if for any $\{x, y\} \subseteq S$ either $x_1 = y_1$ or $x_2 = y_2$ then $S \in L \cup R$.

## 3. The elements of $L^*$

**Proposition 1.** *Let* $A = [x \in \{xy, zx\}]$ *and* $B = [xyz \in \{x, z\}]$. *Then* $A = B = L \cup R$.

*Proof* Let $S \in A$. For any $\{x, y\} \subseteq S$, $x \in \{x(yx), (xy)x\}$ and so $x = xyx$. Therefore, $S \in RB$. Also, $xy \in \{(xy)x, y(xy)\} = \{x, y\}$ and so, by Result 2, $S \in L \cup R$. So $L \cup R \subseteq A \subseteq L \cup R \subseteq B$. Suppose that $T \in B$. For any $\{x, y\} \subseteq T$, $\{(xy)^3, (xy)^4\} \subseteq \{xy\}$ and so $xy = (xy)^2 = x(yx)y \in \{x, y\}$. Then, since $B \subseteq [xyx = x] = RB$, by Result 2, $B \subseteq L \cup R$. ∎

**Proposition 2.** *Let* $C = [xy \in \{y, zw\}]$. *Then* $C = Z \cup R$. *Dually,* $[xy \in \{x, zw\}] = L \cup Z$.

*Proof* Clearly, $Z \cup R \subseteq C$. Suppose $S \in C \setminus (Z \cup R)$. Then there exists $\{a,b,c,d\} \subseteq S$ such that $ab \neq cd$ *and* $cd \neq d$. But $cd \in \{d, ab\}$, which is a contradiction. So $C \subseteq Z \cup R \subseteq C$. ∎

**Proposition 3.** *Let* $V = [xyz \in \{x, z, uv\}$. *Then* $V = L \cup Z \cup R$.

*Proof* Clearly $L \cup Z \cup R \subseteq V$. Suppose $S \in V \setminus (L \cup R)$. Then, by proposition 1, there exists $\{a,b,c\} \subseteq S$ such that $abc \notin \{a,c\}$. But since $S \in V$, for any $uv \in S$, $abc = uv$ and so using Result 1, $S \in Z$. ∎

**Proposition 4.** *Let* $V = [xyz \in \{z, xw\}]$. *Then* $V = \mathbf{I}L \cup R$. *Dually,* $[xyz \in \{x, wz\}] = L \cup \mathbf{I}R$.

*Proof* Clearly, $\mathbf{I}L \cup R \subseteq V$. Let $\{a,b,c,d\} \subseteq S \in V$. Then, $ab(cd) \in \{cd, ab\}$. So, by proposition 1, $S^2 \in L \cup R$. Suppose $S^2 \in R$. If there is some $c \in S - S^2$ then for any $\{a,b\} \subseteq S$, $ac^2 = ac^2c^2 = c^2 \in \{c, ab\}$. Thus, $c^2 = ab$ and so $S^2 = \{0\}$ and $S \in \mathbf{I}L$. If there is no element $c \in S \setminus S^2$, then $S = S^2 \in R$.

Now if $S^2 \in L$ then for any $\{a,b,c\} \in S$, $a^2b^2 = a(ab)b = aab^2 \in \{b, ab\} \cap \{b^2, ab\}$. So either $a^2b^2 = ab$ or $a^2b^2 = b = b^2$. The latter implies that $ab = abb \in \{b, a(ab^2)\} = \{a^2b^2\}$. By Result 1, $S \in \mathbf{I}L$. ∎

**Proposition 5.** *Let* $A = [xyz \in \{xw, qz\}]$, $B = [xy \in \{xuv, zwy\}]$ *and* $D = [xyz = xz ; xyzwuv \in \{xy, uv\}]$. *Then* $A = B = D = \mathbf{I}L \cup \mathbf{I}R$.

*Proof* Let $\{q,u,v,w,x,y,z\} \subseteq S \in A$. Then $(xy)^2 = x(yx)y \in \{xw, qy\}$. Choosing $w = y$ and $q = x$ gives $xy = (xy)^2 \in \{xw, qy\}$ and since $w$ and $q$ are arbitrary elements of $S$, $S \in [xy \in \{xuv, zwy\}] = B$, so $A \subseteq B$. Let $\{u,v,w,x,y,z\} \subseteq S \in B$. Then $xz \in \{xuz, xwz\}$ and so $S \in [xyz = xz = xzxz]$. Since $S \in B$, $(xy)(zw) \in \{(xy)xy, zw(zw)\} = \{xy, zw\}$. By Result 2, $S^2 \in L \cup R$ and so $S \in [xyzwuv \in \{xy, uv\}]$. Hence, $S \in D$ and $B \subseteq D$. If $S \in D$ and $\{x, y\} \subseteq S$ then $xy = (xy)^2$ and, by Result 1 and Proposition 1, $S$ is an inflation of $S^2 \in L \cup R$. So $D \subseteq \mathbf{I}L \cup \mathbf{I}R$. By Result 1, $\mathbf{I}L = [xyz = xw]$ and $\mathbf{I}R = [xyz = qz]$ and so clearly $\mathbf{I}L \cup \mathbf{I}R \subseteq A$. We have shown that $A \subseteq B \subseteq D \subseteq \mathbf{I}L \cup \mathbf{I}R \subseteq A$. ∎

**Proposition 6.** *Let* $V = [xyz \in \{z, xy\}; xyzwuv \in \{xy, uv\}]$. *Then* $V = \mathbf{G}L \cup R$.
*Dually,* $[xyz \in \{x, yz\}; xyzwuv \in \{xy, uv\}] = L \cup \mathbf{G}R$.

*Proof* Using Result 1, clearly $\mathbf{G}L \cup R \subseteq V$. Let $S \in V$. Since $S \in [xyzwuv \in \{xy, uv\}]$, $S^2 \in [xyz \in \{x, z\}]$ and, by proposition 1, $S^2 \in L \cup R$. If $S^2 \in L$ then, by Result 1, $S \in \mathbf{G}L$. So we can assume that $S^2 \in R$. If there exists an element $z \in S \setminus S^2$ then for any $x \in S$, $xzz = xzzxzz = (xzzx)z^2 = z^2$. Since $S \in V$, $xzz = z^2 \in \{z, xz\}$ and, since $z \in S \setminus S^2$, $z^2 = xz (x \in S)$. Let $\{w, q\} \subseteq S$ and choose $x = wq$. Then $z^2 = wqz \in \{z, wq\}$. Therefore, $z^2 = wq$ for any $wq \in S^2$. So $S^2 = \{0\}$. By Result 1, $S \in \mathbf{I}L \subseteq \mathbf{G}L$. If $S \setminus S^2 = \emptyset$ then $S = S^2 \in R$. So we have shown that $V \subseteq \mathbf{G}L \cup R \subseteq V$. ∎

**Proposition 7.** *Let* $V = [xyz \in \{xy, qz\} ; xy \in \{xyz, uvy\}]$. *Then* $V = \mathbf{G}L \cup \mathbf{I}R$.
*Dually,* $[xyz \in \{yz, xq\} ; yz \in \{xyz, yuv\}] = \mathbf{I}L \cup \mathbf{G}R$.

*Proof* By Result 1, clearly $\mathbf{G}L \cup \mathbf{I}R \subseteq V$. Let $S \in V$. Then since $S \in [xy \in \{xyz, uvy\}]$, $xy = (xy)^2$. Since $S \in [xyz \in \{xy, qz\}]$, $xyzw \in \{xy, zwzw\} = \{xy, zw\}$. Therefore $S^2 \in [xy \in \{x, y\}]$. By Result 2, $S^2$ is a chain $C$ of semigroups $S_\alpha (\alpha \in C)$, where $S_\alpha \in L \cup R (\alpha \in C)$. But $S \in [xy \in \{xyz, uvy\}]$ implies that $C$ has only one element. Therefore $S^2 \in L \cup R$. If $S \notin \mathbf{G}L \cup \mathbf{I}R$ then we can assume that $S^2 \in R$, or else $S^2 \in L$ and, by Result 1, $S \in \mathbf{G}L$. Using Result 1 again, since $S \notin \mathbf{I}R$ there exists $\{d,e,f,w\} \subseteq S$ such that $def \neq wf$. But $S \in [xy \in \{xyz, uvy\}]$ and so $wf = wfz (z \in S)$. By substituting $z = def$ and using $S^2 \in R$ we get $wf = wf(def) = def$, a contradiction. So $V \subseteq \mathbf{G}L \cup \mathbf{I}R \subseteq V$. ∎

**Proposition 8.** *Let $A = [xyzwuv \in \{xy, uv\}]$ and $B = [xy \in \{xyw, qxy\}]$. Then $A = B = \mathbf{GL} \cup \mathbf{GR}$.*

*Proof* Let $S \in A$. Then, by proposition 1, $S^2 \in L \cup R$ and, by Result 1, $S \in \mathbf{GL} \cup \mathbf{GR}$. By Result 1, $\mathbf{GL} \cup \mathbf{GR} = \{S : S^2 \in L \cup R\} \subseteq A$. So $A = \mathbf{GL} \cup \mathbf{GR}$. If $S \in B$ then clearly $S^2 \in [x \in \{xy, zx\}] = L \cup R$. By Result 1, $S \in \mathbf{GL} \cup \mathbf{GR}$. So $B \subseteq \mathbf{GL} \cup \mathbf{GR}$. If $S \in \mathbf{GL} \cup \mathbf{GR}$ then for any $\{x, y\} \subseteq S$, $xy = xyxy$. Then, for any $\{q, w\} \subseteq S$, $\{xyw, qxy\} = \{xywxyw, qxyqxy\}$. Now since $S^2 \in L \cup R$, $xy \in \{xywxyw, qxyqxy\} = \{xyw, qxy\}$ So $S \in B$. We have shown that $A = B = \mathbf{GL} \cup \mathbf{GR}$. ∎

**Proposition 9.** *Let $V = [xyx \in \{x, zw\}]$. Then $V = Z \cup RB$.*

*Proof* Clearly $Z \cup RB \subseteq V$. Suppose that $S \in V \setminus RB$. Then there exists $\{a,b\} \subseteq S$ such that $aba \neq a$. Therefore, $aba = zw$ for any $\{z,w\} \subseteq S$. So $S \in Z$. Hence, $V \subseteq Z \cup RB \subseteq V$. ∎

**Proposition 10.** *Let $V = [xyx \in \{x, wx\}]$. Then $V = RB \cup \mathbf{IR}$. Dually, $[xyx \in \{x, xw\}] = \mathbf{IL} \cup RB$.*

*Proof* Let $S \in V$. We show that $S \in [xy = xyxy]$. Let $\{x, y\} \subseteq S$. If $xyx = x$ then $xyxy = xy$. So we can assume that $xyx = yx$. Then $(xy)^2 = (xyx)y = yxy \in \{y, xy\}$. So we can assume that $yxy = y$. But then $(xy)^2 = x(yxy) = xy$, as required. So $S \in [xy = xyxy]$. Now $xyzwxy \in \{xy, (xy)xy\} = \{xy\}$ and by Result 1, $S^2 \in RB$.

Suppose that $T \in V \setminus RB$. Then there exists $\{a,b\} \subseteq T$ such that $aba \neq a$. Therefore $aba = wa (w \in T)$. But then for every $\{u, v, z, q\} \subseteq T$, $zqqa = uvqa$ and so $(uv)_1 = (zq)_1$. This implies that $T^2 \in R$. Let $\{x, y\} \subseteq T$. Since $xy = (xy)^3$, $(xy)_1 = (xyx)_1 \in \{x_1, (x^2)_1\}$. If $x \notin T^2$ then $(xy)_1 = (x^2)_1 = (x^2y^2)_1$. If $x \in T^2$ then $(xy)_1 = (xxy)_1 = (x^2y^2)_1$. So $(xy) = (x^2y^2)_1$. Also, similarly, $yxy \in \{y, y^2\}$ implies $(xy)_2 = (x^2y^2)_2$. Hence $T \in [xy = x^2y^2]$. So $T \in \mathbf{IR}$. Clearly then, $V \subseteq RB \cup \mathbf{IR} \subseteq V$. ∎

**Proposition 11.** *Let $V = \left[ xyx \in \{x, wx, xq\} \; ; \; xy = (xy)^2 \right]$. Then $V = \mathbf{IL} \cup RB \cup \mathbf{IR}$..*

*Proof* Let $\{x, y, z, w\} \subseteq S \in V$. Since $xyzwxy \in \{xy, (xy)^2, (xy)^2\} = \{xy\}$, $S^2 \in RB$. Assume that $S \notin RB$. Then there exists $\{a,b\} \subseteq S$ such that $aba \neq a$. Thus, $aba \in \{wa, aq\}$ for every set $\{w, q\} \subseteq S$. Therefore, $aba \in \{q^*a, aq^*\} \cap \{q^*a, aw^*\} \cap \{w^*a, aw^*\} \cap \{w^*a, aq^*\}$ for every $\{w^*, q^*\} \subseteq S^2$. This gives 16 possible cases and in each case either $(w^*)_1 = (q^*)_1$ or $(w^*)_2 = (q^*)_2$. So by Result 3, $S^2 \in L \cup R$.

Now for any $\{x, y\} \subseteq S$, $x^2y^2x^2 = x(xy^2x)x \in \{x, (xy)x, x(yx)\}$ and so $(x^2y^2)_1 \in \{x_1, (xyx)_1\} = \{x_1, (xy)_1\}$. Therefore, $(x^2y^2)_1 = (xy)_1$. Similarly, $(x^2y^2)_2 = (y^2x^2y^2)_2 \in \{y_2, (yxy)_2\} = \{(xy)_2\}$. So $xy = x^2y^2$ and $S \in \mathbf{IL} \cup \mathbf{IR}$. Clearly then, $V \subseteq \mathbf{IL} \cup RB \cup \mathbf{IR} \subseteq V$. ∎

**Proposition 12.** *Let $V = \left[ xyx \in \{x, xyw\} \; ; \; xy = xyzwxy \right]$. Then $V = \mathbf{GL} \cup RB$. Dually, $\left[ xyx \in \{x, wyx\} \; ; \; xy = xyzwxy \right] = RB \cup \mathbf{GR}$.*

*Proof* Using Result 1, clearly $\mathbf{GL} \cup \mathbf{RB} \subseteq V$. Suppose that $S \in V \setminus \mathbf{RB}$. Then there exists $\{a,b\} \subseteq S$ such that $aba = abw\,(w \in S^2)$. Therefore, $w_2 = q_2$ for every $\{w,q\} \subseteq S^2$. So $S^2 \in L$ and $S \in \mathbf{GL}$. We have shown that $V \subseteq \mathbf{GL} \cup \mathbf{RB} \subseteq V$. ∎

**Proposition 13.** *Let* $V = \left[\, xyx \in \{x, wx, xyq\} \;;\; xy = xyzwxy \,\right]$. *Then* $V = \mathbf{GL} \cup \mathbf{RB} \cup \mathbf{IR}$.
*Dually,* $\left[\, xyx \in \{x, xw, qyx\} \;;\; xy = xyzwxy \,\right] = \mathbf{IL} \cup \mathbf{RB} \cup \mathbf{GR}$.

*Proof* Using Result 1, clearly $\mathbf{GL} \cup \mathbf{RB} \cup \mathbf{IR} \subseteq V$. Let $S \in V \setminus (\mathbf{RB} \cup \mathbf{IR})$. Then, by proposition 10, there exists $\{a,b,w\} \subseteq S$ such that $aba \notin \{a, wa\}$. Therefore, $aba = abq\,(q \in S^2)$. This implies that $u_2 = v_2\,(\{u,v\} \subseteq S^2)$. Therefore, $S^2 \in L$ and $S \in \mathbf{GL}$. We have shown that $\mathbf{GL} \cup \mathbf{RB} \cup \mathbf{IR} \subseteq V \subseteq \mathbf{GL} \cup \mathbf{RB} \cup \mathbf{IR}$. ∎

**Proposition 14.** *Let* $V = \left[\, xyx \in \{x, xyw, qyx\} \;;\; xy = xyzwxy \,\right]$. *Then* $V = \mathbf{GL} \cup \mathbf{RB} \cup \mathbf{GR}$.

*Proof* By Result 1, clearly $\mathbf{GL} \cup \mathbf{RB} \cup \mathbf{IR} \subseteq V$. Suppose that $S \in V \setminus \mathbf{RB}$. Then there exists $\{a,b\} \subseteq S$ such that $aba \neq a$. Therefore $aba \in \{abw, qba\} \cap \{abw, wba\} \cap \{abq, wba\} \cap \{abq, qba\}$, for any $\{w,q\} \subseteq S^2$. This yields 16 cases and in each case either $w_1 = q_1$ or $w_2 = q_2$. So by Result 3, $S^2 \in L \cup R$. By Result 1, $S \in \mathbf{GL} \cup \mathbf{GR}$. ∎

**Proposition 15.** *Let* $V = \left[\, xy \in \{xwy, xyq\} \,\right]$. *Then* $V = \mathbf{GL} \cup \mathbf{IRB}$. *Dually,* $\left[\, xy \in \{xwy, qxy\} \,\right] = \mathbf{IRB} \cup \mathbf{GR}$.

*Proof* By Result 1, clearly $\mathbf{GL} \cup \mathbf{IRB} \subseteq V$. Suppose that $\{x,y,z,w\} \subseteq S \in V \setminus \mathbf{IRB}$. Since $xy \in \{x(yzwx)y, xy(zwxy)\} = \{xyzwxy\}$, by Result 1, $S^2 \in \mathbf{RB}$. Now Result 1 implies there exists $\{a,b,c\} \subseteq S$ such that $acb \neq ab$ and so $ab = abq\,(q \in S^2)$. Therefore for any $w \in S^2$, $q_2 = w_2$ and so $S^2 \in L$ and $S \in \mathbf{GL}$. Hence, $\mathbf{GL} \cup \mathbf{IRB} \subseteq V \subseteq \mathbf{GL} \cup \mathbf{IRB}$. ∎

**Proposition 16.** *Let* $V = \left[\, xy \in \{xyw, xvy, qxy\} \,\right]$. *Then* $V = \mathbf{GL} \cup \mathbf{IRB} \cup \mathbf{GR}$.

By Result 1, clearly $\mathbf{GL} \cup \mathbf{IRB} \cup \mathbf{GR} \subseteq V$. Suppose that $S \in V \setminus \mathbf{IRB}$. Then there exists $\{a,b,c\} \subseteq S$ such that $acb \neq ab$. Therefore, for any $\{w,z\} \subseteq S^2$, $ab \in \{abw, wab\} \cap \{abw, zab\} \cap \{abz, wab\} \cap \{abz, zab\}$. In the 16 possible cases, either $(ab)_1 = w_1 = z_1$ or $(ab)_2 = w_2 = z_2$. So, by Result 3, $S^2 \in L \cup R$. Then, by Result 1, $S \in \mathbf{GL} \cup \mathbf{GR}$. We have shown that $\mathbf{GL} \cup \mathbf{IRB} \cup \mathbf{GR} \subseteq V \subseteq \mathbf{GL} \cup \mathbf{IRB} \cup \mathbf{GR}$. ∎

**Proposition 17.** *Let* $V = \left[\, xyz \in \{yz, xywz\} \;;\; xyzwxy = xy \,\right]$. *Then* $V = \mathbf{GRB}_l \cup \mathbf{GR}$.
*Dually,* $\left[\, xyz \in \{xy, xwyz\} \;;\; xyzwxy = xy \,\right] = \mathbf{GL} \cup \mathbf{GRB}_r$.

*Proof* By Result 1, clearly $\mathbf{GRB}_l \cup \mathbf{GR} \subseteq V$. Suppose that $S \in V \setminus (\mathbf{GRB}_l \cup \mathbf{GR})$. Then by Result 1, $S^2 \in \mathbf{RB}$ and there exists $\{a,b,c,d,e,f,g\} \subseteq S$ such that $abc \neq bc$ and $def \neq degf$. But then, since $a(bcde)f \in \{(bcde)f, a(bcde)gf\}$, either $(abc)_1 = (bc)_1$ (which is a contradiction, as it implies that $abc = bc$) or $(ef)_2 = (gf)_2$, which is also a contradiction as it implies that $def = degf$.
So we have shown that $\mathbf{GRB}_l \cup \mathbf{GR} \subseteq V \subseteq \mathbf{GRB}_l \cup \mathbf{GR}$. ∎

**Proposition 18.** *Let* $V = \left[ xyz \in \{xywz, xqyz\} \; ; \; xy = (xy)^2 \right]$. *Then* $V = \mathbf{GRB}_l \cup \mathbf{GRB}_r$.

*Proof* By definition $\mathbf{GRB}_l \cup \mathbf{GRB}_r \subseteq V$. Suppose that $S \in V \setminus (\mathbf{GRB}_l \cup \mathbf{GRB}_r)$. Then for any $\{l, m, n, x, y, z\} \subseteq S$, $xyz \in \{xy(zlmnxy)z, x(yzlmnx)yz\} = \{xyzlmnxyz\}$ and so $S^3 \in RB$. But $xy = (xy)^2$ and so $S^2 = S^3 \in RB$. Now there exists $\{a,b,c,d,e,f,g,h\} \subseteq S$ such that $abc \neq abdc$ and $efg \neq ehfg$. Now $(efg)(abc) = e(fgab)c \in \{e(fgab)(abd)c, e(hfg)(fgab)c\}$. This implies that either $(abc)_2 = (abdc)_2$ or $(efg)_1 = (ehfg)_1$. But since $(abc)_1 = (abdc)_1$ and $(efg)_2 = (ehfg)_2$, this is a contradiction. Hence, $\mathbf{GRB}_l \cup \mathbf{GRB}_r \subseteq V \subseteq \mathbf{GRB}_l \cup \mathbf{GRB}_r$. ∎

## 4. The lattice of inclusion classes below *I*

Propositions 1-18 give the figure shown in Diagram 1 on page 13. To show that this is the lattice of inclusion classes below *I* we must prove that there is no proper inclusion class between any two points in the diagram. This we proceed to do. In the Propositions below $C = [L_k \subseteq M_k]_{k=1}^t$ is an arbitrary inclusion class.

**Proposition 19.** *If* $a \in S \setminus S^2$ *and* $S \in C$ *then for any* $k \in \{1, 2, ..., t\}$, $x \in L_k \Rightarrow x \in M_k$.

*Proof* Suppose that $L = x \in L_k$ for some $k \in \{1, 2, ..., t\}$ and that $x \notin M_k$. Let $* = (x \to a \; ; \; y \to a^2 \text{ if } y \neq x)$. Now $a = L^* = W_L^* \in M_k$. If $W_L \neq x$ then $a = W_L^* = a^n$ for some $n \geq 2$. This is a contradiction because $a \in S \setminus S^2$ and so $W_L = x \in M_k$. ∎

**Corollary A.** *If* $S \in (Z \setminus \{0\}) \cap C$ *then for any* $k \in \{1, 2, ..., t\}$, $x \in L_k \Rightarrow x \in M_k$.

**Proposition 20.** *If* $S \in (Z \setminus \{0\}) \cap C$ *then* $Z \subseteq C$.

*Proof* Since $S \in (Z \setminus \{0\})$ there exists an element $a \in S \setminus S^2$. If there exists a word $L = x_1 x_2 ... x_n \in L_k$ $(n \geq 2)$ then $W_L$ must have length $\geq 2$; otherwise $0 = L^{*a} = W_L^{*a} = a$, a contradiction. This fact, along with Corollary A implies that any $T \in Z$ satisfies $L_k^* \subseteq M_k^* \; (k \in \{1, 2, ..., t\})$ for any substitution $*$. Hence, $Z \subseteq C$. ∎

**Proposition 21.** *Suppose that* $S \in (L \setminus \{0\}) \cap C$. *Let* $k \in \{1, 2, ..., t\}$. *If* $x_1 x_2 ... x_n \in L_k$ $(n \geq 1)$ *then there exists* $x_1 y_2 y_3 ... y_m \in M_k$ $(m \geq 1)$

*Proof* Suppose $S \in (L \setminus \{0\}) \cap C$ and $L = x_1 x_2 ... x_n \in L_k$ $(n \geq 1)$. Assume that every word $y_1 y_2 ... y_m \in M_k$ $(m \geq 1)$ Satisfies $y_1 \neq x_1$. Let $\{a, b\} \subseteq S$, with $a \neq b$. Let $* = (x_1 \to a \; ; \; y \to b \text{ when } y \neq x_1)$. Then $a = L^* = W_L^* = b$, a contradiction. Hence, $W_L$ must begin in the variable $x_1$. ∎

**Proposition 22.** *If* $S \in (L \setminus \{0\}) \cap C$ *then* $L \subseteq C$.

*Proof* This follows immediately from Proposition 21. ∎

**Proposition 23.** *If* $S \in (R \setminus \{0\}) \cap C$ *then* $R \subseteq C$.

*Proof* Proposition 23 is the dual of Proposition 22. ∎

**Proposition 24.** *If* $S \in (\mathbf{IL} \setminus L) \cap C$ *then* $x \in L_k$ *implies* $x \in M_k$ $(k \in \{1, 2, ..., t\})$.
*Dually, if* $S \in (\mathbf{IR} \setminus R) \cap C$ *then* $x \in L_k$ *implies* $x \in M_k$ $(k \in \{1, 2, ..., t\})$.

*Proof* If $S \in (\mathbf{IL} \setminus L) \cap C$ then there exists $a \in S \setminus S^2$. The result then follows from Proposition 19. ∎

**Proposition 25.** *If* $S \in [\mathbf{IL} \setminus (L \cup Z)] \cap C$ *then* $\mathbf{IL} \subseteq C$. *Dually, if* $S \in [\mathbf{IR} \setminus (R \cup Z)] \cap C$ *then* $\mathbf{IR} \subseteq C$.

*Proof* By Result 1, since $S \in \mathbf{IL} = [xyz = xw]$ and $S \notin L$, there exists $\{a,b\} \subseteq S$ such that $a^2 = ab \neq a \notin S^2$. Suppose that $L = x_1 x_2 ... x_n \in L_k$ $(n \geq 2)$.

Let $c \in S^2$ and suppose that $* = (x_1 \to a\,; y \to c, \text{ if } y \neq x_1)$. Since $S \in \mathbf{IL} = [xyz = xw]$, $ac = a^2$ and $L^* = (x_1 x_n)^* = a^2$. Assume that $W_L$ is *not* of the form $x_1 y_2 y_3 ... y_m$ $(m \geq 2)$. Then $W_L^* \subseteq \{a, c\}$. But since $a \neq a^2$, $a^2 = L^* = W_L^* \subseteq \{a, c\}$ and so $c = a^2$ $(c \in S^2)$. This implies that $S \in Z$, a contradiction. Hence $W_L$ is of the form $x_1 y_2 y_3 ... y_m$ $(m \geq 2)$ and $W_L \in M_k$. Together with Proposition 24 this proves that any $T \in \mathbf{IL} = [xyz = xw]$ is a member of $C$. ∎

**Proposition 26.** *Suppose that* $Z \subseteq C$. *If* $S \in (\mathbf{GL} \setminus \mathbf{IL}) \cap C$ *then* $\mathbf{GL} \subseteq C$.
*Dually, if* $Z \subseteq C$ *and* $S \in (\mathbf{GR} \setminus \mathbf{IR}) \cap C$ *then* $\mathbf{GR} \subseteq C$.

*Proof* Suppose that $L = x_1 x_2 ... x_n \in L_k$ $(k \in \{1, 2, ..., t\})$. Note that since $S \in \mathbf{GL} \setminus \mathbf{IL} = [xyz = xy] \setminus [xyz = xw]$, there exists $\{a,b,c\} \subseteq S$ such that $ab \neq ac$, with $a \in S \setminus S^2$. Suppose that there is no word $W \in M_k$ of the form $W = x_1 x_2 y_3 ... y_m$ $(m \geq 2)$. Define $* = (x_1 \to a\,; x_2 \to b\,; y \to a^2, \text{ if } y \notin \{x_1, x_2\})$ and $\bullet = (x_1 \to a\,; x_2 \to c\,; y \to a^2, \text{ if } y \notin \{x_1, x_2\})$. Then $ab = L^* = W_L^* \in M_k^* = \{a,b,a^2,b^2,ba\}$ and $ac = L^\bullet = W_L^\bullet \subseteq M_k^\bullet = \{a,c,a^2,c^2,ca\}$.

Now we show that $ab \in \{a,b,a^2,b^2,ba\}$ implies that $ab = a^2$. Recall that $S \in [xyz = xy] = \{T : T^2 \in L\}$. Firstly, since $a \in S \setminus S^2$, $ab \neq a$. If $ab = b$ then $ab = a(ab) = a^2b = a^2$. If $ab = b^2$ then $ab = abb = ab^2 = a(ab) = a^2b = a^2$. Finally, if $ab = ba$ then $ab = aba = aab = a^2$. So $ab = a^2$. Similarly, $ac = a^2$. But then $ab = ac$, a contradiction. Hence, $M_k$ contains a word of the form $W = x_1 x_2 y_3 ... y_m$ $(m \geq 2)$. Since $Z \subseteq C$, proposition 19 implies that any $T \in \mathbf{GL} = [xyz = xy]$ is a member of $C$. ∎

**Proposition 27.** *If* $S \in [RB \setminus (L \cup R)] \cap C$ *then* $RB \subseteq C$.

*Proof* Note that since $S \in [RB \setminus (L \cup R)] \cap C$, by Result 1 there exists $\{a,b\} \subseteq S$ such that $ab \notin \{a,b\}$. Note also that $S \in RB \subseteq [xyz = xz]$. Let $k \in \{1, 2, ..., t\}$. Suppose that $x = L \in L_k$. Then $x \in M_k$, or else with $* = (x \to a\,; y \to b, \text{ if } y \neq x)$, $L^* = x^* = a = W_L^* \in \{b, ab, ba\}$, which contradicts the fact that $ab \notin \{a,b\}$.

Now assume $L = x_1 x_2 ... x_n \in L_k (n \geq 2)$, with $x_1 = x_n$. Let $* = (x_1 \to a \,;\, y \to b$, if $y \neq x_1)$. Then $L^* = (x_1^2)^* = a = W_L^*$. Therefore, $W_L$ must begin and end in the variable $x_1$, or else $a \in \{ab, ba\}$, which contradicts $ab \notin \{a, b\}$. So $W_L = x_1 y_2 y_3 ... y_{m-1} x_1 \, (m \geq 1)$.

Now assume $L = x_1 x_2 ... x_n \in L_k (n \geq 2)$. Assume that there is no word of the form $x_1 y_2 y_3 ... y_{m-2} y_{m-1} x_n \in M_k (m \geq 2)$. Then let $* = (x_1 \to ab \,;\, x_n \to ba \,;\, y \to b$, if $y \notin \{x_1, x_n\})$. So $abba = a = L^* = W_L^* \in M_k^* \subseteq \{b, ab, ba\}$, contradicting the fact that $ab \notin \{a, b\}$. Therefore, there exists a word $x_1 y_2 y_3 ... y_{m-2} y_{m-1} x_n \in M_k (m \geq 2)$.

This implies that any $T \in RB$ is a member of $C$. ∎

**Proposition 28.** *If* $S \in [IRB \setminus (IL \cup RB \cup IR)] \cap C$ *then* $IRB \subseteq C$.

*Proof* Since $S \in IRB \setminus RB$ there exists $c \in S \setminus S^2$. By Proposition 11, there exists $\{a, b, w, q\} \subseteq S \in [xyz = xz]$ such that $a^2 = aba \notin \{a, wa, aq\}$ and $a \notin S^2$.

Let $L = x \in L_k \, (k \in \{1, 2, ..., t\})$. Then, by Proposition 19, $x \in M_k$.

Suppose $n > 1$ and $L = x_1 x_2 ... x_n \in L_k \, (k \in \{1, 2, ..., t\})$. Assume there is no word $x_1 y_2 y_3 ... y_{m-2} y_{m-1} x_n \in M_k (m \geq 2)$. Let $* = (x_1 \to a \,;\, x_n \to a \,;\, y \to wq$, if $y \notin \{x_1, x_n\})$. Recall that $S \in [xyz = xz]$. Then $a^2 = L^* = W_L^* \in M_k^* \subseteq \{a, wq, wa, aq\}$. So $a^2 = wq = w^2 q^2$.

Let $* = (x_1 \to a \,;\, x_n \to a \,;\, y \to qw$, if $y \notin \{x_1, x_n\})$. Then $a^2 = L^* = W_L^* \in M_k^* \subseteq \{a, qw, aw, qa\}$. So $a^2 \in \{qw, aw, qa\}$. If $a^2 = qw$ then $a^2 = (a^2)^2 = (w^2 q^2)(qw) = w^2 = a^3 = w^2 a = wa$, a contradiction. If $a^2 = aw$ Then $a^2 = (a^2)^2 = (w^2 q^2)(aw) = w^2 = a^3 = w^2 a = wa$, a contradiction. If $a^2 = qa$ then $a^2 = (a^2)^2 = (qa)(w^2 q^2) = q^2 = a^3 = aq^2 = aq$, a contradiction. Therefore, there exists a word $x_1 y_2 y_3 ... y_{m-2} y_{m-1} x_n \in M_k (m \geq 2)$. Hence, any $T \in IRB$ is a member of $C$. ∎

**Proposition 29.** *If* $S \in [GRB_l \cap C] \setminus [GL \cup IRB]$ *then* $GRB_l \subseteq C$.
*Dually,* $S \in [GRB_r \cap C] \setminus [IRB \cup GR]$ *implies* $GRB_r \subseteq C$.

*Proof* Since $S \notin GL \cup IRB$, by proposition 15 there exists $\{a, b, c, d, e\} \subseteq S$ such that $(ab)_1 \neq (a^2)_1$.

Note that $S \in GRB_l = [xyz = xywz \,;\, xy = (xy)^2]$, which implies by Result 1 that $S^2 \in RB$. Note that for any $\{x, y\} \subseteq S$, $xyy = xy(x)y = (xy)^2 = xy$ and so $(xy)_2 = (y^2)_2$. Therefore $(y^2)_2 = (xy)_2 = (x^2 y^2)_2$. This implies that there exists $\{a, b\} \subseteq S$ such that $(ab)_1 \neq (a^2)_1$, or else by Result 1, $S \in [xy = x^2 y^2] = IB$ and so $S \in IRB$, a contradiction. Also, $(ab)_1 \neq (a^2)_1$ implies $a \notin S^2$ and so Proposition 19 is valid.

Furthermore, since $S \notin GL$ there exists $\{c,d,e\} \subseteq S$ such that $ec \neq ecd$. Clearly $(ec)_1 = (ecd)_1$. Therefore $(c^2)_2 \neq (d^2)_2$.

Suppose that $L = x^2 \in L_k$ ($k \in \{1, 2, ..., t\}$). Let $* = (x \to a; y \to bc^2$, if $y \neq x)$. Then $a^2 = L^* = W_L^*$ and so $a^2 = a^3 = aW_L^*$. Therefore $W_L$ begins in the variable $x$, or else $(ab)_1 = (a^2)_1$. Also, $W_L$ must have length $\geq 2$, or else $a^2 = L^* = W_L^* = a$, a contradiction. Then, if $W_L = w_1 w_2 ... w_n$ ($n \geq 2$) then $w_1 = x = w_2$, or else $(ab)_1 = (a^2)_1$. Thus, $W_L \in \{x^2, x^2 Kx, x^2 Ky\}$ ($y \neq x$), where $K$ may be the empty word. Note that if $W_L = x^2 Ky$ then $a^2 = L^* = W_L^* = a^2 K^* bc^2$ and then $(c^2)_2 = (a^2)_2$. Similarly, either $M_k$ contains a word of the form $x^2$ or $x^2 Kx$ or $(d^2)_2 = (a^2)_2$. So we have shown that either $(c^2)_2 = (a^2)_2 = (d^2)_2$, which is a contradiction, or $M_k$ contains a word of the form $x^2$ or $x^2 Kx$.

Now suppose that $L = xy \in L_k$ ($k \in \{1, 2, ..., t\}$) and $x \neq y$. Let $* = (x \to a; y \to bd; z \to a^2 c$, if $z \notin \{x, y\})$. Then $L^* = (xy)^* = abd = W_L^*$, where $W_L = w_1 w_2 ... w_n \in M_k$. If $w_1 = y$ then $abd = W_L^* = bdW$ for some (possibly empty) word $W$. Then $ab = abab = ab(dW)b = a(bdW)b = a(abd)b$ and so $(ab)_1 = (a^2)_1$, a contradiction. Clearly, if $w_1 = z \notin \{x, y\}$ then $(ab)_1 = (a^2)_1$. So $w_1 = x$. Also, $n \geq 2$ or else $a \in S^2$.

Suppose that $W_L$ has length 2. Then if $w_2 \in \{x, z\}$ ($z \notin \{x, y\}$) it easily follows that $(ab)_1 = (a^2)_1$. So $W_L = xy$. Assume that $n \geq 3$. Then since $S \in [xyz = xywz]$, $W_L = w_1 w_2 w_n$. It is then straightforward to show that $w_1 = x$ and $w_2 = y$. If $w_n \in \{x, z\}$ ($z \notin \{x, y\}$) then $L^* = (xy)^* = abd = W_L^* \in \{abda, abda^2 c\}$ and so $(d^2)_2 \in \{(a^2)_2, (c^2)_2\}$. Hence, $(d^2)_2 = (a^2)_2$.

Similarly, using $\bullet = (x \to a; y \to bc; z \to a^2 d$, if $z \notin \{x, y\})$ we can show that either $(c^2)_2 = (a^2)_2$ or $W_L = xyKy$, for some (possibly empty) word $K$. So either $(c^2)_2 = (a^2)_2 = (d^2)_2$, a contradiction, or $M_k$ contains a word $W_{L,*} \in \{xy, xyKy\}$ for some substitution $*$.

Now assume that $L = x_1 x_2 ... x_n \in L_k$ ($k \in \{1, 2, ..., t\}$), where $n \geq 3$. Then since $S \in GRB_l = [xyz = xywz]$, $L = x_1 x_2 x_n$. Define $o = |\{x_1, x_2, x_n\}|$, the order of the set $\{x_1, x_2, x_n\}$. We consider the three cases where $o \in \{1, 2, 3\}$.

Case 1. ($o = 1$) Then define $x = x_1 (= x_2 = x_n)$ and $* = (x \to a; y \to bc^2$, if $y \neq x)$. Then $L^* = a^3 = a^2 = W_L^*$ for some word $W_L \in M_k$, with $W = w_1 w_2 ... w_m$ ($m \geq 1$). If $m = 1$ then $w_1 \neq x$ implies $a^2 = W_L^* = bc^2$ and so $a^2 = a^3 = a(bc^2)$ and $(ab)_1 = (a^2)_1$, a contradiction. If $w_1 = x$ then $a^2 = W_L^* = a$, a contradiction also. So $m \geq 2$.

Now $w_1 = w_2 = x$, or else $(a^2)_1 = (ab)_1$ and if $w_m \neq x$ then $(a^2)_2 = (c^2)_2$. So we have shown that $M_k$ contains a word $W_L \in \{x^2, x^2 Kx\}$, for some possibly empty word $K$, or $(a^2)_2 = (c^2)_2$. Similarly, using the substitution $* = (x \to a \, ; \, y \to bd^2$, if $y \neq x)$ we can show that $M_k$ contains a word $W_L \in \{x^2, x^2Kx\}$, for some possibly empty word $K$, or $(a^2)_2 = (d^2)_2$. Since $(d^2)_2 \neq (c^2)_2$, $M_k$ contains a word $W_L \in \{x^2, x^2Kx\}$, for some possibly empty word $K$.

Case 2. ($o = 2$)

Case 2.1 ($x = x_1 = x_2 \neq x_n = y$) Define $* = (x \to a \, ; \, y \to bc^2 \, ; \, z \to bd^2$, if $z \notin \{x, y\})$. Now $L^* = a^2bc^2 = a^2c = W_L^*$ for some word $W_L = w_1w_2...w_m \in M_k$ $(m \geq 1)$.

Suppose $m = 1$. If $w_1 = x$ then $a^2c = a \in S^2$, a contradiction. If $w_1 = y$ or $w_1 = z \notin \{x, y\}$ then $(a^2)_1 = (ab)_1$, also a contradiction. So $m \geq 2$. If $m \geq 2$ then $w_1 = x = w_2$, or else $(a^2)_1 = (ab)_1$. Hence $m = 2$ implies $W_L = x^2$ and so $(a^2)_2 = (c^2)_2$. Similarly, if $\bullet = (x \to a \, ; \, y \to bd^2 \, ; \, z \to bc^2$, if $z \notin \{x, y\})$ then we can prove that $(a^2)_2 = (d^2)_2$. Then $(d^2)_2 = (c^2)_2$, a contradiction, and so $m \neq 2$. So we can assume that $m \geq 3$.

We have already proved that $w_1 = x = w_2$. If $w_m = x$ then $a^2c = W_L^* = (xx...x)^* = a^3 = a^2$ and so $(c^2)_2 = (a^2)_2$. If $w_m = z \notin \{x, y\}$ then $a^2c = W_L^* = (xx...z)^* = a^2...bd^2$, which implies $(c^2)_2 = (d^2)_2$, a contradiction. Hence, either $(c^2)_2 = (a^2)_2$ or $M_k$ contains a word $W_L = xxKy$, for some (possibly empty) word $K$. Similarly, using the substitution $\bullet = (x \to a \, ; \, y \to bd^2 \, ; \, z \to bc^2$, if $z \notin \{x, y\})$ we can prove that either $(d^2)_2 = (a^2)_2$ or $M_k$ contains a word $W_L = xxKy$, for some (possibly empty) word $K$. Since $(d^2)_2 \neq (c^2)_2$, we conclude that $M_k$ contains a word $W_L = xxKy$, for some (possibly empty) word $K$.

Case 2.2 ($x = x_1 \neq x_2 = x_n = y$). Then $L = xyy = xy$ and this case has already been dealt with in the earlier part of this proof.

Case 2.3 ($x = x_1 = x_n \neq x_2 = y$). Then $L = xyx$. Define $* = (x \to a \, ; \, y \to bc \, ; \, z \to a^2c$, if $z \notin \{x, y\})$. Then $L^* = aba = W_L^*$, for some word $W_L = w_1w_2...w_m$ $(m \geq 1)$. If $m = 1$ then $w_1 = x$, or else $(a^2)_1 = (ab)_1$, a contradiction. But then $L^* = aba = W_L^* = a \in S^2$, a contradiction. So $m \neq 1$.

Suppose that $m = 2$. Then $W_L = xy$, or else $(a^2)_1 = (ab)_1$, a contradiction. Hence, $aba = abc$ and so $(a^2)_2 = (c^2)_2$.

Suppose that $m \geq 3$. Then $w_1 = x$ and $w_2 = y$, or else $(a^2)_1 = (ab)_1$, a contradiction. Also, either $w_m = x$ or $(a^2)_2 = (c^2)_2$. So either $(a^2)_2 = (c^2)_2$ or $M_k$ contains a word of the form $W_{L,*} = xyKx$, for some (possibly empty) word $K$.

Similarly, using the substitution $\bullet = (x \to a \, ; \, y \to bd \, ; \, z \to a^2d$, if $z \notin \{x, y\})$ we can prove that either $(a^2)_2 = (d^2)_2$ or $M_k$ contains a word of the form $W_{L,\bullet} = xyKx$, for some (possibly empty) word $K$. Then since $(c^2)_2 \neq (d^2)_2$, $M_k$ contains a word of the form $W = xyKx$, for some (possibly empty) word $K$.

Case 3 $(o = 3)$ $(x = x_1 \neq x_2 = y \neq x_n = z \neq x)$. Then $L = xyz$.

Define $* = (x \to a\,;\, y \to bd\,;\, z \to a^2c\,;\, w \to a^2d$, if $w \notin \{x, y, z\})$. Then $L^* = abda^2c = abc = W_L^*$ for some word $W_L = w_1 w_2 ... w_m \in M_k$ $(m \geq 1)$. If $m = 1$ then $w_1 = x$, or else $(a^2)_1 = (ab)_1$, a contradiction. But then $a = abc \in S^2$, a contradiction. So $m \neq 1$.

If $m = 2$ then $W_L = xy$, or else $(a^2)_1 = (ab)_1$, a contradiction. But then $abc = W_L^* = abd$ and $(c^2)_2 = (d^2)_2$, also a contradiction. So $m \geq 3$.

Now $w_m \notin \{y, w\}$ $(w \notin \{x, y, z\})$, or else $(c^2)_2 = (d^2)_2$, a contradiction. If $w_m = x$ then $(a^2)_2 = (c^2)_2$. So either $(a^2)_2 = (c^2)_2$ or $W_{L,*} = xyKz$, for some (possibly empty) word $K$. Similarly, using the substitution

$\bullet = (x \to a\,;\, y \to bc\,;\, z \to a^2d\,;\, w \to a^2c$, if $w \notin \{x, y, z\})$ we can prove that either $(a^2)_2 = (d^2)_2$ or $W_{L,\bullet} = xyKz$, for some (possibly empty) word $K$. Since $(c^2)_2 \neq (d^2)_2$, there exists a word $W = xyKz \in M_k$, for some (possibly empty) word $K$.

So we have proved the following: for any $k \in \{1, 2, ..., t\}$,

1) $x \in L_k \Rightarrow x \in M_k$ ; 2) $x^2 \in L_k \Rightarrow M_k$ contains a word of the form $x^2$ or $x^2Kx$;

3) $xy \in L_k \,(x \neq y) \Rightarrow M_k$ contains a word of the form $xy$ or $xyKy$ ;

4) $xxHx \in L_k \Rightarrow M_k$ contains a word of the form $x^2$ or $x^2Kx$;

5) $xxHy \in L_k \,(x \neq y) \Rightarrow M_k$ contains a word of the form $x^2Ky$ ;

6) $xyHx \in L_k \,(x \neq y) \Rightarrow M_k$ contains a word of the form $xyKx$;

7) $xyHy \in L_k \,(x \neq y) \Rightarrow M_k$ contains a word of the form $xyKy$ and

8) $xyHz \in L_k \,(x \neq y \neq z \neq x) \Rightarrow M_k$ contains a word of the form $xyKz$.

Thus, for any $T \in \mathbf{GRB}_l = \left[ xyz = xywz\,;\, xy = (xy)^2 \right]$ we have proved that for any word $L \in L_k \,(k \in \{1, 2, ..., t\})$ there exists a word $W \in M_k$ such that for any transformation $*$, $L^* = W^*$. Hence, $\mathbf{GRB}_l \subseteq C$. ∎

## 5. Proof of the main Theorem

We now proceed to prove that Diagram 1 below is a diagram of the lattice of inclusion classes under the inclusion class $\mathbf{I} = \left[ xyz \in \{xywz, xqyz\}\,;\, xy = (xy)^2 \right]$. To do this we must prove that there is no inclusion class between any two points in the diagram, other than those indicated in the diagram. First we list some definitions.

**Definition.** If $\alpha = (A_1, A_2, ..., A_n, A_{n+1})$ is a sequence of $n+1$ distinct points in a lattice then $\alpha$ is an **n-path** from $A_1$ to $A_{n+1}$ if $A_1 < A_2 < ... < A_n < A_{n+1}$.

**Definition.** If $\alpha = (A, B)$ is a sequence of two distinct points in a lattice **L** then **B covers A** if $\alpha$ is a 1-path from $A$ to $B$ and for any $X \in \mathbf{L}$, $A < X < B$ implies $X \in \{A, B\}$. If **B** covers 0 then **B** is called an atom.

**Theorem 30.** *For any $n \geq 1$, if $A_1 < X < A_{n+1}$ and $\alpha = (A_1, A_2, ..., A_n, A_{n+1})$ is an n-path from $A_1$ to $A_{n+1}$ in Diagram 1 then $X$ appears as a point in Diagram 1. If $n = 1$ then $A_2$ covers $A_1$.*

Proof. Examining Diagram 1 we see that, using Propositions 20, 22, 23, 25, 26, 27, 28 and 29, the Theorem holds for $n = 1$.

The proof for $n \geq 2$ also follows from Propositions 20, 22, 23, 25, 26, 27, 28 and 29 but the proof is very long. We will give an idea of the methods used by proving that if $0 < X < (GRB_l \cup GRB_r)$ then $X$ appears as a point in Diagram 1.

DIAGRAM 1.

THE LATTICE OF INCLUSION CLASSES UNDER
$$I = \left[ xyz \in \{xywz, xqyz\}; xy = (xy)^2 \right]$$

(continuation of the proof of Theorem 30.)

Since $0 < X < (GRB_l \cup GRB_r)$ there are four cases.
Case 1. There exists $S \in [X \cap GRB_r] \setminus (IRB \cup GR)$ and there exists $T \in [X \cap GRB_l] \setminus (GL \cup IRB)$;
Case 2. There exists $S \in [X \cap GRB_r] \setminus (IRB \cup GR)$ and there is no $T \in [X \cap GRB_l] \setminus (GL \cup IRB)$;
Case 3. There is no $S \in [X \cap GRB_r] \setminus (IRB \cup GR)$ and there exists $T \in [X \cap GRB_l] \setminus (GL \cup IRB)$;
Case 4. There is no $S \in [X \cap GRB_r] \setminus (IRB \cup GR)$ and there is no $T \in [X \cap GRB_l] \setminus (GL \cup IRB)$.

Case 1. In this case, by Proposition 29, $(\mathbf{GRB}_l \cup \mathbf{GRB}_r) \subseteq X$ and so $X = (\mathbf{GRB}_l \cup \mathbf{GRB}_r)$.

Case 2. In this case, by Proposition 29, $\mathbf{GRB}_r \subseteq X \subseteq (\mathbf{GRB}_l \cup \mathbf{GRB}_r)$. Now if $X \cap [\mathbf{GRB}_l \setminus \mathbf{GRB}_r] = \emptyset$
Then $X = \mathbf{GRB}_r$. If $X \cap [\mathbf{GRB}_l \setminus \mathbf{GRB}_r] \neq \emptyset$ then there exists $T \in X \cap (\mathbf{GL} \cup \mathbf{IRB}) \cap [\mathbf{GRB}_l \setminus \mathbf{GRB}_r]$
Then, since $(\mathbf{IRB} \cup \mathbf{IL}) \subseteq \mathbf{GRB}_r$, $T \in X \cap (\mathbf{GL} \setminus \mathbf{IL})$. By Proposition 26, $\mathbf{GL} \subseteq X$ and so
$\mathbf{GL} \cup \mathbf{GRB}_r \subseteq X \subseteq (\mathbf{GRB}_l \cup \mathbf{GRB}_r)$. Since $(\mathbf{GRB}_l \cup \mathbf{GRB}_r)$ covers $\mathbf{GL} \cup \mathbf{GRB}_r$
$X \in \{\mathbf{GL} \cup \mathbf{GRB}_r, (\mathbf{GRB}_l \cup \mathbf{GRB}_r)\}$.

We have proved that $X \subseteq \{\mathbf{GRB}_r, \mathbf{GL} \cup \mathbf{GRB}_r, (\mathbf{GRB}_l \cup \mathbf{GRB}_r)\}$.

Case 3. As in the proof of Case 2, $X \in \{\mathbf{GRB}_l, \mathbf{GRB} \cup \mathbf{GR}, (\mathbf{GRB}_l \cup \mathbf{GRB}_r)\}$.

Case 4. In this case $X \in \mathbf{GL} \cup \mathbf{IRB} \cup \mathbf{GR}$.

We define $X_1 = X \cap [\mathbf{IRB} \setminus (\mathbf{IL} \cup \mathbf{RB} \cup \mathbf{IR})]$, $X_2 = X \cap [\mathbf{RB} \setminus (\mathbf{L} \cup \mathbf{R})]$, $X_3 = X \cap (\mathbf{GL} \setminus \mathbf{IL})$,
$X_4 = X \cap (\mathbf{GR} \setminus \mathbf{IR})$, $X_5 = X \cap [\mathbf{IL} \setminus (\mathbf{L} \cup \mathbf{Z})]$, $X_6 = X \cap [\mathbf{IR} \setminus (\mathbf{R} \cup \mathbf{Z})]$, $X_7 = X \cap (\mathbf{L} \setminus \mathbf{0})$,
$X_8 = X \cap (\mathbf{R} \setminus \mathbf{0})$, $X_9 = X \cap (\mathbf{Z} \setminus \mathbf{0})$ and $X_{10} = \{\mathbf{0}\}$.

Clearly, $X = \bigcup X_i \, (i \in \{1,2,...,9,10\})$. Using Propositions 20, 22, 23, 25, 26, 27, 28 and 29 it follows that
$X_1 \neq \emptyset$ implies $\mathbf{IRB} \subseteq X$, $X_2 \neq \emptyset$ implies $\mathbf{RB} \subseteq X$, $X_3 \neq \emptyset$ implies $\mathbf{GL} \subseteq X$, $X_4 \neq \emptyset$ implies $\mathbf{GR} \subseteq X$,
$X_5 \neq \emptyset$ implies $\mathbf{IL} \subseteq X$, $X_6 \neq \emptyset$ implies $\mathbf{IR} \subseteq X$, $X_7 \neq \emptyset$ implies $\mathbf{L} \subseteq X$, $X_8 \neq \emptyset$ implies $\mathbf{R} \subseteq X$ and
$X_9 \neq \emptyset$ implies $\mathbf{Z} \subseteq X$. This implies that $X$ is one of the inclusion classes that appear in Diagram 1.
*This completes the outline of the proof of Theorem 30.* ∎

**Corollary 31.** *Diagram 1 is the lattice of inclusion classes under the inclusion class $\mathbf{I}$. Its elements are unions of varieties of generalized inflations of rectangular bands.*

Note that not every semigroup inclusion class is a union of varieties of semigroups, as we proceed to prove.

**Proposition 32.** *The semigroup inclusion class $[xy \in \{x,y\}]$ is not a union of semigroup varieties.*

Proof. Let S be the two element chain. If we assume that $[xy \in \{x,y\}]$ is a union of semigroup varieties then $S \in V \cap [xy \in \{x,y\}]$ for some semigroup variety $V \subseteq [xy \in \{x,y\}]$. Let $V = [u_\alpha = v_\alpha]_{\alpha \in J}$. Then for any $\alpha \in J$ the words $u_\alpha$ and $v_\alpha$ must have equal content (that is, the words must consist of the same variables). Consider the set $T = \{0, x, y\}$ with the following multiplication: $x = xx$, $y = yy$, and
$0 = 00 = 0x = x0 = 0y = y0 = xy = yx$. Then T is a semilattice and so $T \in V \setminus [xy \in \{x,y\}]$, a contradiction. ∎

**Open Questions**

1. Does $\mathbf{GRB}$ cover $\mathbf{I}$ ?
2. Describe the modular, upper modular, lower modular, co-distributive and distributive elements of $\mathbf{L}^*$.

[1] R. A. R. Monzo
10 Albert Mansions
Crouch Hill, London N8 9RE
United Kingdom
bobmonzo@talktalk.net